\documentclass{amsart}

\usepackage{amssymb}

\usepackage{url}

\newtheorem{theorem}{Theorem}[section]

\newtheorem{remark}[theorem]{Remark}
\theoremstyle{definition}

\newtheorem{definition}[theorem]{Definition}

\newcommand{\Z}{\mathbb{Z}}

\newcommand{\C}{\mathbb{C}}

\newcommand{\msk}{\medskip}

\newcommand{\nin}{\noindent}

\DeclareMathOperator\arcsinh{arcsinh}
\DeclareMathOperator\arccosh{arccosh}
\begin{document}

%\linenumbers

\title[Sin, cos, exp and log of Liouville Numbers]{Sin, cos, exp, and log of Liouville Numbers}


%    author one information

\author{Taboka Prince Chalebgwa and Sidney A. Morris}
\address{The Fields Institute for Research in Mathematical Sciences,
222 College Street, Toronto, Ontario, MST 3J1,
Canada}
\email{taboka@aims.ac.za}
\address{School of Engineering, IT and Physical Sciences,
Federation University Australia,
PO Box 663,
Ballarat, Victoria, 3353,
Australia
and\newline
Department of Mathematical and Physical Sciences,
La Trobe University,
Melbourne, Victoria, 3086,
Australia}
\email{morris.sidney@gmail.com}

%    \subjclass is required.

\date{}

\dedicatory{\nin\hspace{0.6in}Dedicated to  Kurt Mahler}

\begin{abstract} For any Liouville number $\alpha$, all of the following are transcendental  numbers: 
$\textrm{e}^\alpha$, $\log_\textrm{e}\alpha$, $\sin \alpha$, $\cos\alpha$, $\tan\alpha$, $\sinh\alpha$,  $\cosh\alpha$, $\tanh\alpha$,  $\arcsin\alpha$  and the inverse functions evaluated at $\alpha$ of the listed trigonometric and hyperbolic functions, noting that  wherever multiple values are involved, every such value is transcendental.  This remains true if ``Liouville number'' is replaced by ``$U$-number'', where $U$ is one of  Mahler's classes of transcendental numbers.
\end{abstract}

\maketitle

\section{Introduction}  

\noindent In 1844 Joseph Liouville proved the existence of transcendental numbers \cite{Angell, Baker}. He introduced the set $\mathcal L$ of real numbers, now known as Liouville numbers, and showed that they are all transcendental.

\begin{definition}\label{d2}

Recall that a real number $\xi$ is called a Liouville number if for every positive integer $n$, there exists a pair of integers $(p,q)$ with $q > 1$, such that 

\begin{equation*}
0 < \left| \xi - \frac{p}{q} \right| < \frac{1}{q^n}. \eqno{\qed}
\end{equation*}

\end{definition}

Alan  Baker in his classic  work \cite{Baker} on transcendental number theory   said: A classification of the set of all transcendental number into three disjoint aggregates, termed $S$-, $T$-, and $U$-numbers was introduced by Kurt Mahler \cite{Mahler} in 1932, and it has proved to be of considerable value in the general development of the subject. 

In this paper we demonstrate just how powerful and useful Mahler's classification of transcendental numbers is.

The following beautiful theorem, which is a corollary of the Lindemann-Weierstrass Theorem  appears as Theorem 9.11 in Ivan Niven's book \cite{Niven}. We prove the analogous result with ``algebraic number'' replaced by ``Liouville number''.
\begin{theorem}\label{t1.1}
The following numbers are transcendental:
\begin{enumerate}
\item[(i)] $\textrm{e}^\alpha$, $\sin \alpha$, $\cos\alpha$, $\tan\alpha$, $\sinh\alpha$, $\cosh\alpha$, $\tanh\alpha$;
\item[(ii)] $\log_\textrm{e}\alpha, \arcsin\alpha$ and in general the inverse functions of those listed in (i), \end{enumerate} for any non-zero algebraic number $\alpha$; wherever multiple values are involved, every such value is transcendental.  
\end{theorem}

It is not widely known that $\textrm{e}^\alpha$ and $\log_\textrm{e}\alpha$ are transcendental numbers when $\alpha$ is a Liouville number, though the $\exp$ case is stated explicitly in \cite[p.$\,$98]{Kumar} and the $\log$ case, as pointed out to  the second author by Michel Waldschmidt, is  implicit in \cite[\S3.5]{Bugeaud}. The proof in our paper for $\exp$ is different from that in \cite{Kumar}.

On the other hand, the results for trigonometric and hyperbolic functions do not appear in print. The proofs of all these results for Liouville numbers, and indeed a wider class of numbers, depend on properties of the Mahler classes of transcendental numbers.
 \msk

 \section{Mahler Classes}

We follow the presentation in \cite[\S3]{Bugeaud}. While the definitions and results therein are stated and proved for real numbers, \textit{mutatis mutandis}, they carry over to the case of complex numbers. The classification of Mahler partitions the complex numbers into four sets, characterized by the rate with which a nonzero polynomial with integer coefficients approaches zero when evaluated at a particular number.

Given a polynomial $P(X) \in \C[X]$, recall that the height of $P$, denoted by $H(P)$, is the maximum of the absolute values of the coefficients of $P$. Given a complex number $\xi$, a positive integer $n$, and a real number $H \geq 1$, we define the quantity

\begin{equation*}
w_n(\xi, H) = \min \{ |P(\xi)| : P(X) \in \Z[X], \ H(P) \leq H, \ \deg (P) \leq n, \ P(\xi) \neq 0 \}.
\end{equation*}

Furthermore, we set

\begin{equation*}
w_n(\xi) = \limsup_{H \rightarrow \infty} \frac{- \log w_n(\xi, H)}{\log H}
\end{equation*}

and

\begin{equation*}
w(\xi) = \limsup_{n \rightarrow \infty} \frac{w_n(\xi)}{n}.
\end{equation*}

With the above notation in mind, Mahler partitions the complex numbers into the following:

\begin{definition}\label{d1}

Let $\xi$ be a complex number. We say $\xi$ is an

\begin{itemize}
\item $A$-number if $w(\xi) = 0$,
\item $S$-number if $0 < w(\xi) < \infty$,
\item $T$-number if $w(\xi) = \infty$ and $w_n(\xi) < \infty$ for any $n \geq 1$,
\item $U$-number if $w(\xi) = \infty$ and $w_n(\xi) = \infty$ for all $n \geq n_0$, for some positive integer $n_0$. \qed

\end{itemize}

\end{definition}

\begin{remark}\label{r1}
Note that, the $A$-numbers are the algebraic numbers. \qed
\end{remark}

The following theorem of Mahler, see \cite[Theorem 3.2]{Bugeaud}, records a fundamental property of  the Mahler classes.

\begin{theorem}\label{t2.3}
If $\xi, \eta \in \C$ are algebraically dependent, then they belong to the same Mahler class.

\end{theorem}

The following theorem, which is key to our main result, was proved by Kurt Mahler \cite{Mahler, Mahler1953}.  (See \cite[\S3.5]{Bugeaud}.)

\begin{theorem}\label{t2.4}  If $a$ is an  algebraic number with $a\ne0$, then $\exp(a)\in  S$ and if $a\ne0,1$  then $\log(z)\in S\cup T$.
\end{theorem} 

Mahler also proved in \cite{Mahlerpi} the following result which we shall need.

\begin{theorem}\label{t2.4a} The number $\pi\in S\cup T$.
\end{theorem}

\begin{remark}\label{2.7}
Note that the Liouville numbers are $U$-numbers. Furthermore, if $\xi$ is a Liouville number, then $i \xi$ is a $U$-number by Theorem \ref{t2.3}. \qed
\end{remark}

\section{The Main Theorem}

\begin{theorem}\label{3.1} For any  $U$-number $\alpha$, in particular for $\alpha$ any Liouville number,  all of the following are transcendental numbers: 
$\textrm{e}^\alpha$, $\log_\textrm{e}\alpha$, $\sin \alpha$, $\cos\alpha$, $\tan\alpha$, $\sinh\alpha$,  $\cosh\alpha$, $\tanh\alpha$ and the inverse functions evaluated at $\alpha$ of the listed trigonometric and hyperbolic functions, noting that  wherever multiple values are involved, every such value is transcendental. 
\end{theorem}

\begin{proof}
For ease of notation, we adopt the conventional notation that $\log$ denotes $\log_\textrm{e}$. We shall demonstrate the result for $e^\alpha$, $\sin \alpha$, $\tan \alpha$, $\sinh \alpha$ and their  inverses. The corresponding proofs for the remaining functions in each family are analogous.

Henceforth, let $\alpha$ be  a $U$-number.

\begin{enumerate}

\item $\exp \alpha$: Suppose that $\mu = e^\alpha$ is an algebraic number. Then $\log \mu  = \alpha \in S \cup T$ by Theorem \ref{t2.4}. This is a contradiction since $\alpha \in U$. So $e^\alpha$ is a transcendental number.

\item $\log \alpha$: Suppose $\mu$ is one of the values of $ \log \alpha$ and is an algebraic number. Then $e^\mu  = \alpha \in S$ by Theorem \ref{t2.4}. This is a contradiction since $\alpha \in U$. So all values of $\log \alpha$ are  transcendental.

\item $\sin \alpha$: Suppose $\mu = \sin \alpha$ is an algebraic number. Then by Remark \ref{2.7} and Theorem \ref{t2.4}, we have that $i \alpha \in U$. By (1)  $t = e^{i \alpha}$ is transcendental.  Further,  $2i \sin \alpha = e^{i \alpha} - e^{-i \alpha} = t - \frac{1}{t} = 2i\mu = \beta$, where, by our supposition, $\beta$ is an algebraic number. Now $t - \frac{1}{t} = \beta$  implies that $t = \frac{\beta \pm \sqrt{4 - \beta^2}}{2}$. Since $\beta$ is an algebraic number, the right hand side of the preceding equation is also an algebraic number, and this is a contradiction since $t$ is a transcendental number. So $\sin \alpha$  is transcendental.

\item $\arcsin \alpha$: Suppose firstly that one of the values of $\arcsin \alpha = \mu$ is $0$. Then $\alpha=k\pi$, for some $k\in \mathbb{Z}$.  If $k=0$, then $\alpha=0$ which contradicts the fact  that $\alpha\in U$. If $\alpha=k\pi$, $k\ne0$, then $k\pi\in S\cup T$ by Theorem~\ref{t2.4a}, which contradicts the fact  that $\alpha\in U$.

Next suppose  that one of the values of $\arcsin \alpha = \mu$ is an algebraic number.  
Recall that $\arcsin \alpha = -i \log (i \alpha + \sqrt{1-\alpha^2})$.  
Now $i \mu  = \log (i \alpha + \sqrt{1-\alpha^2})$ implies  $e^{i \mu} = \sqrt{1-\alpha^2} + i \alpha$. By Theorem \ref{t2.4}, $e^{i \mu} \in S$. 
Putting $X = \alpha$ and $Y = \sqrt{1-\alpha^2} + i \alpha$ we have that  they satisfy the equation $P(X,Y) = Y^4 + 4Y^2X^2 - 2Y^2 + 1 = 0$, and hence $X$ and $Y$ are algebraically dependent. So by Theorem \ref{t2.3}, $X$ and $Y$ are in the same Mahler class. As we are given $X\in U$, this implies $e^{i\mu}= Y\in U$ which is a contradiction as we saw that  $e^{i\mu}\in S$. Thus  all values of $\arcsin \alpha$ are transcendental.

\item $\tan \alpha$: Suppose $\mu = \tan \alpha$ is algebraic. Then this implies that $i \mu = \frac{t - 1/t}{t + 1/t}$ is algebraic, where $t = e^{i \alpha}$, which is transcendental by Theorem \ref{t2.4}. The former equation implies that $t^3 + t + i \beta t - i\beta = 0$, where $\beta = \frac{1}{\mu}$. In the interest of brevity, we omit exhibiting the general solution for $t$, but we note that the polynomial has algebraic coefficients, hence, each solution $t$ is also algebraic, a contradiction.

\item $\arctan \alpha$: Put $\mu = \arctan \alpha$, so $2i \mu = \log \left( \frac{i-\alpha}{i + \alpha} \right)$. Note that, as in  (4), $\mu \neq 0$. Suppose $\mu$ is an  algebraic number. Then $e^{2i \mu} \in S$ by Theorem \ref{t2.4}. On the other hand, $X = \alpha$ and $Y = \frac{i - \alpha}{i + \alpha}$ satisfy the equation $P(X,Y) = X^2Y^2 + 2X^2Y + X^2 - Y^2 + 2Y - 1 = 0$, hence are algebraically dependent. Theorem \ref{t2.3} then implies that $\frac{i - \alpha}{i + \alpha} \in U$, a contradiction. So every value of  $\arctan \alpha$ is transcendental.

\item $\sinh \alpha$: Recall that $\sinh \alpha = \frac{e^\alpha - e^{-\alpha}}{2}$. By  (1), $t = e^\alpha$ is transcendental. Suppose $t - \frac{1}{t} = 2 \mu$ is an  algebraic number.  This implies that $t = \mu \pm \sqrt{1 - \mu^2}$ is also algebraic, a contradiction. So  $\sinh \alpha$ is transcendental.

\item $\arcsinh \alpha$: We proceed as in (4). Suppose $\mu = \arcsinh \alpha = \log (\alpha + \sqrt{\alpha^2 + 1})$ is algebraic. By Theorem \ref{t2.4}, we have that $e^\mu \in S$. However $X = \alpha$ and $Y = \alpha + \sqrt{1+ \alpha^2}$ satisfy the equation $P(X,Y) = Y^2 - 2XY -1 = 0$. Hence $X$ and $Y$  are algebraically dependent and therefore $\alpha + \sqrt{1+ \alpha^2} \in U$, a contradiction. So $\arcsinh \alpha$ is transcendental.

\end{enumerate}

\end{proof}

\begin{remark}\label{3.2} In fact  the above argument shows that if $\alpha$ is in the Mahler class   $T$, then $\log(\alpha)$ is a transcendental number. Additionaly, the theorem remains true for the composition of a trigonometric or hyperbolic function with the inverses of the other functions in the corresponding family. For instance, if $\alpha$ is a Liouville number, then $\sinh (\arccosh \alpha)$ is transcendental. \qed
\end{remark}

\begin{remark}\label{3.3}\rm
We conclude by recording that Corollary 6 of \cite{Waldschmidt} implies that $\exp(\alpha)$ is a Louiville number for an uncountable number of Liouville numbers $\alpha$. Recall Edmond T. Maillet's result, \cite[Chapter 3]{Maillet}, which says that if $t$ is a Liouville number  and $R(x)$ is a rational function with rational coefficients, then $R(t)$ is a Liouville number. In our case we use 
$R(t) = \frac{1}{2} (t+\frac{1}{t})$.  
It   follows from this and \cite{Waldschmidt} that there exists an uncountable set of Liouville  numbers $\alpha$ such that $\sinh \alpha$ is a Liouville number and there exists an uncountable set of Liouville numbers $\alpha$ such that $\cosh \alpha$ is a Liouville number. It is not obvious if this statement  remains true when $\sinh$ is  replaced by $\sin$ or $\cos$. \qed\end{remark}

\end{document}